\documentclass[a4paper,11pt]{article}
\usepackage{amsmath}
\usepackage{amssymb}
\usepackage{natbib}
\usepackage{graphicx}
\PassOptionsToPackage{hyphens}{url}
\usepackage{hyperref}
\usepackage{cleveref}
\usepackage{titlesec}
\usepackage{titling}
\crefformat{section}{\S#2#1#3}
\crefformat{subsection}{\S#2#1#3}
\crefformat{subsubsection}{\S#2#1#3}
\titleformat{\paragraph}
{\normalfont\normalsize\bfseries}{\theparagraph}{1em}{}
\titlespacing*{\paragraph}
{0pt}{3.25ex plus 1ex minus .2ex}{1.5ex plus .2ex}

\begin{document}
\begin{center}
\LARGE{On finding all positive integers $a,b$ such that $b\pm a$ and $ab$ are palindromic}
\\
$ $
\\
\large{Wang Pok Lo, Yuval Paz}
\date{}
\end{center}

\begin{abstract}
It is proven that the only integer solutions $(a,b)$ such that $a+b$ and $ab$ are palindromic are $(2,5\cdot 10^k-3)$, $(3,24)$ and $(9,9)$, and in a similar fashion, $b-a$ and $ab$ are only palindromic at $(a,b)=(3,147\cdot 10^{4(k+1)}+5247\sum_{i=0}^k10^{4i})$, $(3,161\,247\cdot 10^{4k+7}+5247\sum_{i=0}^k10^{4i+3}+387)$, $(3,147)$ and $(3,161\,247\,387)$ for $k=0,1,2,\cdots$. Note $a\le b$ without loss of generality.
\end{abstract}

\section{Introduction}

The challenge to determine all positive integers $a,b$ such that $b+a$ and $ab$ are palindromic has been explored by a few people, but none has yet provided a rigorous proof of all the solutions. In 2009, the conjectured solutions were posted on OEIS by Mark Nandor [2]. More recently, there was a question on Quora [3] asking this, and we have even done so on Mathematics Stack Exchange where a user Michael Lugo [4] conjectured the same as Nandor. In this paper, we will prove their claims, and will also generalise this to the case where $b-a$ and $ab$ are palindromic.
\\
\\
\textit{Definition}: 

\textit{Two integers are palindromic if the digits of one integer are the same as the reverse of the digits of the other, and if they both have the same number of digits.}
\\
\\
That is, if one integer $a$ has digit representation $a_na_{n-1}\cdots a_1$ and the other integer $b$ has digit representation $b_nb_{n-1}\cdots b_1$, then they are palindromic if and only if $$a_1=b_n, \quad a_2=a_{n-1},\quad \cdots , \quad a_{n-1}=b_2, \quad a_n = b_1$$
Furthermore, $a_2,\cdots,a_{n-1},b_2,\cdots,b_{n-1}\in\{0,1,2,\cdots,9\}$ but $a_1,a_n,b_1,b_n\in\{1,2,\cdots,9\}$. This is so that the the first and last digits of $a$ and $b$ are not zero; otherwise, they would not have the same number of digits.

\section{When $a+b$ and $ab$ are palindromic}

\subsection{Formulating the problem}\label{sec:2.1}

Without loss of generality assume that $a\le b$. Then $\frac{ab}{a+b} > \frac{ab}{2b}=\frac{a}{2}$ and this must be less than $10$. Hence $a<20$.
\\
\\
For $n>1$, let $b$ have the above digit representation, and let 
\begin{equation}
a+b=c_nc_{n-1}\cdots c_2c_1=10^{n-1}c_n+10^{n-2}c_{n-1}+\cdots+10c_2+c_1
\end{equation}
\begin{equation}
ab=c_1c_2\cdots c_{n-1}c_n=10^{n-1}c_1+10^{n-2}c_2+\cdots+10c_{n-1}+c_n
\end{equation}
Substituting $b$ from (1) into (2), we get
\begin{equation}
(a\cdot 10^{n-1}-1)c_n+(a\cdot 10^{n-2}-10)c_{n-1}+\cdots+(10a-10^{n-2})c_2+(a-10^{n-1})c_1=a^2
\end{equation}
Note that to preserve the same number of digits in each expression, $ac_n<10$, except for extreme cases highlighted in $\mathsection$\ref{sec:2.3}.

\subsection{Finding the solutions}\label{sec:2.2}
In this section it will be assumed that $c_n=b_n$. Exceptions to this are also discussed in $\mathsection$\ref{sec:2.3}.
\subsubsection{Solutions when $a=2$}

If $n=1$, we can immediately solve the equation $2+b=2b\implies b=2$ so $(2,2)$ is a trivial solution.
\\
\\
From (3), the equation becomes
\begin{equation}
(2\cdot 10^{n-1}-1)c_n+(2\cdot 10^{n-2}-10)c_{n-1}+\cdots+(20-10^{n-2})c_2+(2-10^{n-1})c_1=4
\end{equation}
Clearly $c_n$ is even, and since $2\cdot 6>10$, $c_n = 2,4$. If $c_n$ is the former, then $b_n=2$ so the first digit of $2b$ is $c_1=4$ or $5$. If $c_1=4$, then $b_1=c_1-2=2$ so the last digit of $2b$ must be $4=c_n$ which is a contradiction. Similarly, if $c_1=5$, then $b_1=3$ so the last digit of $2b$ must be $6=c_n$, again, a contradiction. Therefore $c_n=4$, meaning that $2b$ starts with either $c_1=8,9\implies b_1=6,7$. If $b_1=6$ then the last digit of $2b$ is $2$ contradicting the fact that the first digit of $b$ is $4$. Hence $b_1=7\implies c_1=9$.
\\
\\
Dividing equation (4) by $2$, the RHS is still even, so to fulfill that on the LHS, we must have that $-5c_{n-1}+c_1$ is even since all other terms on that side have at least one even factor. Since $c_1$ is odd, so is $c_{n-1}$. Notice that $c_1=2c_n+1$, implying that there is carrying. This narrows $c_{n-1}$ down to either being $5$, $7$ or $9$.
\\
\\
If $c_{n-1}=5$, then $2b$ ends in the digits $54$ so $b$ must end in the digits $27$ or $77$, implying that $c_2=2,7$. However the second digit of $2b$ is $c_2=2c_{n-1}=0,1$ which is a contradiction. Similarly, if $c_{n-1}=7$, then $2b$ ends in the digits $74$ so $b$ must end in the digits $37$ or $87$, implying that $c_2=3,8$. However the second digit of $2b$ is $c_2=2c_{n-1}=4,5$ which is a contradiction. Finally, if $c_{n-1}=9$, then $2b$ ends in the digits $94$ so $b$ must end in the digits $47$ or $97$, implying that $c_2=4,9$. However the second digit of $2b$ is $c_2=2c_{n-1}=8,9$ which forces $c_2=9$.
\\
\\
We have now arrived at $b=49\cdots 97\implies 2+b=49\cdots99$ and $2b=99\cdots 94$ so $c_3\ge 5$. Again, we previously showed that $c_1$ is odd so $c_{n-2}$ is also odd. This is a cycle, so the only solutions when $a=2$ are $b=2,47,497,4997,\cdots$ which can be generalised to $5\cdot 10^k-3$ for $k=0,1,2,\cdots$.

\subsubsection{Solutions when $a=3$}

In this section, equation (3) will be used for $n>2$. That said, the cases $n=1,2$ will firstly be considered. Of course, $3+b=3b$ gives no integer solutions so this eliminates the first one. If $n=2$, equation (3) can be modified to give $29c_2-7c_1=9$ which is a standard Diophantine equation. Solving using the Euclidean Algorithm gives the general solution $(c_1,c_2)=(9+7t,36+29t)$ for an integer $t$. But since $c_1,c_2<10$, the only possible solution is when $t=-1$, so $(c_1,c_2)=(2,7)\implies b=24$.
\\
\\
For $n>2$, we have
\begin{equation*}
1+(3\cdot 10^{n-1}-1)c_n+(3\cdot 10^{n-2}-10)c_{n-1}+\cdots+(30-10^{n-2})c_2+(3-10^{n-1})c_1=10
\end{equation*}
so $1-c_n+3c_1$ is a multiple of $10$. As $3c_n<10$, $c_n$ is restricted to $1,2,3$. 
\\
\\
If $c_n=1$, $10$ divides $1-1+3c_1$ which is impossible. 
\\
\\
If $c_n=2$, $10$ divides $3c_1-1$. This can be achieved only if $c_1=7$, so $3b$ starts with $7$. Since the first digit of $3+b$ is $2$, this indicates carrying, and in particular, that $c_{n-1}=3,4,5,6$. If $c_{n-1}=3$, $3b$ ends in $32$ so $b$ ends in $44$. This means that $b+3$ ends in $47$ and in turn, $3b$ ends in $74$. A contradiction arises as $3\cdot 23<74$. If $c_{n-1}=4$, $3b$ ends in $42$ so $b$ ends in $14$. This means that $b+3$ ends in $17$ and in turn, $3b$ starts with $71$. However this is impossible as $3\cdot 24>71$. If $c_{n-1}=5$, $3b$ ends in $52$ so $b$ ends in $84$. This means that $b+3$ ends in $87$ and in turn, $3b$ starts with $78$. Again this is contradictory since $3\cdot 25<78$. Finally, if $c_{n-1}=6$, $3b$ ends in $62$ so $b$ ends in $54$. This means that $b+3$ ends in $57$ and in turn, $3b$ starts with $75$ which is impossible as $3\cdot 26>75$. No solutions exist in this category.
\\
\\
If $c_n=3$, $10$ divides $3c_1-2$. This can be achieved only if $c_1=4$. Now $3c_n=9=c_1$ only since carrying will increase the number of digits so immediately there is a contradiction. Therefore the only solution when $a=3$ is $b=24$.

\subsubsection{Solutions when $a=4,5$}

If $n=1$ then it is easy to show that no solutions exist for $a=4,5$.
\\
\\
For $a=4$, $c_n=1,2$ and we have that
\begin{equation*}
(4\cdot 10^{n-1}-1)c_n+(4\cdot 10^{n-2}-10)c_{n-1}+\cdots+(40-10^{n-2})c_2+(4-10^{n-1})c_1=16
\end{equation*}
and clearly $c_n$ must be even to keep the parities consistent on both sides of the equation. Thus if $c_n=2$, $c_1=8,9$ including carrying leading to $b_1=4,5$. However, $4b_1$ ends in $6,0$ respectively, not $c_n=2$, which is a contradiction.
\\
\\
Similarly, for $a=5$, $c_n=1$ and we have that
\begin{equation*}
(5\cdot 10^{n-1}-1)+(5\cdot 10^{n-2}-10)c_{n-1}+\cdots+(50-10^{n-2})c_2+(5-10^{n-1})c_1=25
\end{equation*}
but the LHS is not divisible by $5$. This is again a contradiction.

\subsubsection{Solutions when $a=7,9$}

When $a=7$, $c_n=1$ so that $ac_n<10$ so equation (3) becomes
\begin{equation*}
(7\cdot 10^{n-1})+(7\cdot 10^{n-2}-10)c_{n-1}+\cdots+(70-10^{n-2})c_2+(7-10^{n-1})c_1=50
\end{equation*}
Now for $n>2$, every term on the LHS is divisible by $10$ except $7c_1$ so $c_1$ must be divisible by $10$. This is a contradiction as $c_1$ must be a single digit and cannot be zero. The case where $n=1$ trivially gives no solutions and if $n=2$, this requires the solution of the Diophantine equation $(7\cdot 10-1)c_2-(7-10)c_1=7^2\implies 3(23c_2-c_1)=49$, but the RHS is not divisible by $3$. Hence no solutions exist and this completes $a=7$.
\\
\\
Similarly, when $a=9$, $c_n=1$ so
\begin{equation*}
(9\cdot 10^{n-1})+(9\cdot 10^{n-2}-10)c_{n-1}+\cdots+(90-10^{n-2})c_2+(9-10^{n-1})c_1=82
\end{equation*}
For $n>2$, every term on the LHS is even except $9c_1$ so $c_1$ must be even. This means that $c_n$ is also even which is a contradiction. The case where $n=1$ trivially also gives no solutions and this case where $n=2$ requires the solution of the Diophantine equation $(9\cdot 10-1)c_2+(9-10)c_1=9^2\implies 89c_2-c_1=81$, which is $(c_1,c_2)=(t,-81+89t)$. As $0<c_1,c_2<10$, the solution occurs when $t=1$, so $(c_1,c_2)=(1,8)$, implying that $b=18-9=9$. Hence the only solution when $a=9$ is $b=9$.

\subsubsection{Solutions for the remaining $a$}

If $a$ is even, so is the RHS of (3). The LHS is also even as $a$ and $10$ are divisible by $2$, except for the term $-c_n$. This means that $c_n$ is even. However, for even $a>5$, we must set $c_n=1$ so that $ac_n<10$ which is again a contradiction. If $a$ is odd; that is, $a=11,13,\cdots,19$, then clearly $ac_n>10$ so all that is left is to check the extreme cases in the next section, where, for example, $b=10^n-a$ so that $a+b$ and $ab$ have the same number of digits.

\subsection{Checking the extreme cases}\label{sec:2.3}

\subsubsection{When $a=2,3,4,5,7,9$}

There were instances in $\mathsection$\ref{sec:2.2} where it was assumed that $c_n=b_n$. While this may be true for the majority of the values of $b$, there are still some exceptions. For example, consider the case $a=2$. If $b=19\cdots 98$ or $19\cdots 99$ then $c_n$ becomes $2$ not $1$. However, by inspection it is apparent that $2+b$ and $2b$ are not palindromic for these $b$, since $2+b$ starts with $2$ but $2b$ ends in $6$ or $8$. Similarly, for the other five values of $a$, we can check from $b=10^{n-1}(b_n+1)-a$ to $10^{n-1}(b_n+1)-1$ - since here, $c_n=b_n+1$, but this yields no solutions either for all $n>1$.

\subsubsection{When $a=11,13,15,17,19$}

For these values of $a$, $b_n=9$ so that $a+b$ and $ab$ have the same number of digits. This means that we need only check from $b=10^n-a$ to $10^n-1$ and it can be easily verified (by comparing the first and last digits of each of $a+b$ and $ab$) that no solutions exist here either.

\section{When $b-a$ and $ab$ are palindromic}
\subsection{Formulating the problem}

This is very similar to section 2.1. The LHS of equation (1) in $\mathsection$\ref{sec:2.1} will be replaced by $b-a$ but the rest of (1) and (2) remain the same. We get a near equivalent equation to (3); the only difference is due to the negative sign on the RHS:
\begin{equation}
(a\cdot 10^{n-1}-1)c_n+(a\cdot 10^{n-2}-10)c_{n-1}+\cdots+(10a-10^{n-2})c_2+(a-10^{n-1})c_1=-a^2
\end{equation}
Note that $a<10$; otherwise, $ab$ will have at least one more digit than $b-a$ as subtraction of a positive integer cannot increase the value of the expression. Again, the criterion (bar exceptions) that $ac_n<10$ still holds.

\subsection{Finding the solutions}

As in $\mathsection$\ref{sec:2.2}, assume that $c_n=b_n$ in this section.

\subsubsection{Solutions when $a=2,4,5$}
We will start with $a=5$. We have two cases: $c_n=0$ and $c_n=5$, as $5b$ ends in $c_n$. Of course, $c_n= 0$ is trivially false, and $c_n=5$ implies that $5b$ has more digits than $b-5$. Both cases lead to contradictions, so there are no solutions.
\newline
\newline
We will now consider $a=4$, first assuming $n=1$ we get $b-4=4b$, and this implies no solutions. Now we can apply (5) and get:
\begin{equation*}
(4\cdot 10^{n-1}-1)c_n+(4\cdot 10^{n-2}-10)c_{n-1}+\cdots+(40-10^{n-2})c_2+(4-10^{n-1})c_1=-16
\end{equation*}
Therefore $c_n$ is even, and to keep $ac_n<10$ we also get $c_n<3\implies c_n=2$. Multiplying by $4$ we get $c_1=8,9$ (which is the last digit of $b-4$) and we need only consider the carrying of $1$ because beyond that we will get a new digit. Now $b=2\cdots\cdots2$ or $b=2\cdots\cdots3$, and multiplying the possible last digits by $4$ we get $c_n=8$ for the former and $c_n=2$ in the latter. Thus, if solutions exists, $c_n=2$ and $c_1=3$.
\\
\\
This gives $b=2\cdots\cdots3\implies b-4=2\cdots\cdots9\implies 4b=9\cdots\cdots2$; in other words, there is a carrying of $1$ from $c_{n-1}$, so we get two possible values: $c_{n-1}=3,4$.
\newline
If $c_nc_{n-1}=23$ then $4b$ ends with $32$, dividing $32$ by $4$ yields $b_1=8$, which contradicts the fact that $b_1=3$.
\newline
If $c_nc_{n-1}=24$ then $4b$ ends with $42$ which is not divisible by $4$ so we have another contradiction. So for $a=4$ there are no solutions.
\newline
\\
For $a=2$, $n=1$ gives the equation $2a=a-2$ which yields no results, and if $n\ge 2$, we can start by noticing that $c_n=2,4$ because
\begin{equation*}
(2\cdot 10^{n-1}-1)c_n+(2\cdot 10^{n-2}-10)c_{n-1}+\cdots+(20-10^{n-2})c_2+(2-10^{n-1})c_1=4
\end{equation*}
If $c_n=4$, $c_1$ is either $8$ or $9$, but $c_1$ cannot be $8$ as $8+2$ ends in $0$. Thus $c_1=9\implies b_1=1$, but $2b$ ends in $4$, which is a contradiction.
\\
\\
If $c_n=2$, $c_1$ is $4$ or $5$, but we know that $2b$ ends in $c_n=2$. Hence $b$ ends in either $1$ or $6$, so $b-2$ ends with $9$ or $4$ which is $c_1$. Now the only common number between $4,5$ and $4,9$ is $4$, so $c_1=4$. It is easy to see that $24+2=26$ is not a solution. So the possible values of $c_{n-1}$ are $0,1,2,3,4$.
\\
\\
If $c_{n-1}=0,2,4$ then $2b$ ends with $d2$ where $d=0,2,4$. This is a contradiction as for a positive integer $k$ we have $2\cdot (k\cdot 10+6)$ ends in $t2$ where $t$ is odd.
\newline
\\
If $c_{n-1}=1$ then $b-2=21\cdots\cdots4\implies 2b=4\cdots\cdots12$, which means that $b$ ends in $06$ or $56$, so $c_2=0,5$ which is a contradiction to the fact that $b$ starts with $21$. We are just left with $c_{n-1}=3$; if $b-2=23\cdots\cdots 4$ then $2b=4\cdots\cdots 32$, and $c_{n-1}=3$ also means that $c_2=6,7$. Hence $b=23\cdots\cdots 66$ or $b=23\cdots\cdots 76$, but $76\cdot 2=152$ does not end in $32$ so $c_{2}=6$.
\\
\\
\textit{Claim:}

\textit{Let} $3_k=\underbrace{3\cdots 3}_{k~\text{times}},6_k=\underbrace{6\cdots 6}_{k~\text{times}}$. \text{Then} $b=2(3_k)\cdots6_k$ \text{implies} that $b=2(3_{k+1})\cdots6_{k+1}$.
\\
\\
\textit{Proof:}
\\
\\
To prove this claim it will be shown that $b=2(3_k)6_k$ is not a solution, and indeed there exists the digit $7$ in $2b$ but not in $b-2$.
\\
\\
If $23_k6_k$ is not a solution, there exists more digits, namely $c_{n-k-1},c_{k+1}$, that are not in the $3_k$ or $6_k$, so we need to check $c_{n-k-1}=0,1,2,3,4$.
\\
\\
For even $c_{n-k-1}$ we get an easy contradiction just like at the start: for an even digit $d$,\newline $b-2=2(3_k)d\cdots 6_{k-1}4\implies 2b=4(6_{k-1})\cdots d3_k2$, but this is impossible as the last $6$ has a carrying which makes the next digit odd, and $d$ is even.
\\
\\
For $c_{n-k-1}=1$ we get $b-2=2(3_k)1\cdots 6_{k-1}4\implies 2b=4(6_{k-1})\cdots 13_k2$ that means $b$ ends in $06_{k-1}$ or $56_{k-1}$ and just like we did in the case of $c_{n-1}=1$, it is a contradiction.
\\
\\
This forces $c_{n-k-1}=3$, so $b-2=2(3_{k+1})\cdots 6_{k-1}4\implies2b=4(6_{k-1})\cdots 3_{k+1}2$, and $c_{n-k-1}=3\implies c_{k+1}=6,7$. If $c_{k+1}=7$ we get $5$ in the place where there should be $3$, so $c_{k+1}=3$. We are done, because $6\ne 3$ and the number of digits is ever growing (infinite), but every integer is finite, so there are no solutions. $\square$

\subsubsection{Solutions when $a=6,7,8,9$}

From (5), we get that when $a=6,8$, $c_n$ is even, but if $c_n>1$ then $ac_n>10$, so for $a=6,8$ there are no solutions.
\\
\\
For $a=7$ we get $c_n=1$, and $c_1=7,8,9$. If $c_1=7$ then $b=1\cdots\cdots4$ so $7b=7\cdots\cdots8$, but $8\ne 1$ is a contradiction. For $c_1=8$ we get $b=1\cdots\cdots 5$, then $7b$ ends in 5 and not $1$, again, a contradiction. Lastly, if $c_1=9$ then $b=1\cdots\cdots 6$, so $7b$ ends in $2\ne1$. So there are no solutions. For $a=9$ we have $c_n=1$ and $c_1=9$, so $b=1\cdots\cdots 8$ and $9b$ ends in $2$ and not in $9$, so again no solutions.

\subsubsection{Solutions when $a=3$}

For $n\le 3$, it is possible to solve the respective Diophantine equations, and no solutions exist when $n=1,2$. When $n=3$, the Diophantine equation becomes
\begin{equation*}
(3\cdot 10^2-1)c_3+(30-10)c_2+(3-100)c_1=-3^2\implies 299c_3+20c_2-97c_1=-9
\end{equation*}
and the general solutions are $(c_1,c_2,c_3)=(t,-135-299s+259t,9+20s-17t)$ for integers $s,t$. Since $c_1=t<10$, there are only nine cases to cover. For $t=1,2,3,7,8,9$, $c_3\ge10$ and for $t=5,6$, $c_2\gg10$. However, when $t=4$, $c_2=4$ and $c_3=1$ so $b-3=144\implies b=147$ is the only solution for $n\le 3$.
\\
\\
For the rest of this section, $n>3$. Firstly, we know that $c_n=1,2,3$. Plugging $a=3$ and subtracting one from both sides of (5), the equation becomes
\begin{equation}
-1+(3\cdot 10^{n-1}-1)c_n+(3\cdot 10^{n-2}-10)c_{n-1}+\cdots+(30-10^{n-2})c_2+(3-10^{n-1})c_1=-10
\end{equation}
and it can be seen that $10$ divides $-1-c_n+3c_1$. However, if $c_n=2$ then $3(c_1-1)$ is divisible by $10$ which is impossible as $c_1$ is a single-digit positive integer. If $c_n=3$ this forces $c_1=8$, so $b-3$ (and hence $b$) begins with $3$ and $3b$ begins with $8$ due to their palindromicity. This is contradictory since $3\cdot 3=9>8$, so the only case left is $c_n=1$.
\\
\\
If $c_n=1$ this forces $c_1=4$ since $10$ must divide $3c_1-2$. Now $b$ starts with $1$ and $3b$ starts with $4$ so there is a carrying of $1$ from the second digit to the first digit - implying that $c_{n-1}=4,5,6$. If $c_{n-1}=5$, $b-3=15\cdots\cdots 4\implies b=15\cdots\cdots 7$ and $3b=4\cdots\cdots 51$. This forces $c_2=1$ since $17\cdot3=51$ so $3b$ starts with $41$. We reach a contradiction as $15\cdot3>41$. This leaves us with $c_{n-1}=4,6$.

\paragraph{When $c_{n-1}=4$}

If $c_{n-1}=4$, $b=14\cdots\cdots 47\implies b-3=14\cdots\cdots 44\implies 3b=44\cdots\cdots 41$ and since there is a carrying of $2$ from the third digit, $c_{n-2}=7,8,9$. If $c_{n-2}=9$ then $b_3=6=c_3$ since $647\cdot 3=1941$. Thus $b$ starts with $149$ and $3b$ with $446$ which is a contradiction as $149\cdot 3>446$. Similarly, if $c_{n-2}=8$ then $b=3=9=c_3$ since $947\cdot 3=2841$. Thus $b$ starts with $148$ and $3b$ with $449$ which is again a contradiction as $148\cdot3=444<449$.
\\
\\
If $c_{n-2}=7$, $b_3=2=c_3$ since $247\cdot 3=741$. Thus $b$ starts with $147$ and $3b$ with $442$ so there is a carrying of $1$ from the fourth digit. This implies that $c_{n-3}=4,5,6$. 
\\
\\
If $c_{n-3}=4$, $b=1474\cdots  247\implies 3b=442\cdots 4741$ which forces $c_4=8$ as $8247\cdot 3=24741$. Hence $3b$ starts with $4428$ which is contradictory since $1474\cdot 3=4422$. If $c_{n-3}=6$, $b=1476\cdots 247\implies 3b=442\cdots 6741$ which forces $c_4=2$ as $2247\cdot 3=6741$. Hence $3b$ starts with $4422$ which is contradictory yet again since $1476\cdot 3=4428$. In a similar argument, it can be shown that $c_{n-3}=5=c_4\,(*)$. In the next five paragraphs we will demonstrate that a pattern shows up.
\\
\\
Suppose that $b$ has seven digits; that is, $b=147x247$ where $x$ is a digit from $0$ to $9$. If this is multiplied by $3$, $3b=441(3x)741$. But $3b$ starts with $442$ due to palindromicity, so $3x=12,15,18$ or that $x=4,5,6$. This forces $x=5$ since the last digit of $3x$ is the same as $x$ itself. Hence $b=1475247$ is a solution.
\\
\\
Suppose that $b$ has eight digits; that is, $b=147xy247$ where $x,y$ are single digits. From similar reasoning to the above, we deduce that $3(xy)=yx+100$ (concatenation, not multiplication of $xy$) due to carrying, so $29x-7y=100$. The general solution to this is $(x,y)=(100+7t,400+29t)$ for an integer $t$, but it is impossible for $y$ to be positive and less than $10$. No solutions exist.
\\
\\
Suppose that $b$ has nine digits; that is, $b=147xyz247$. The Diophantine equation this time is $299x+20y-97z=100$ with general solution $(x,y,z)=(20+20s-17t,-249-299s+259t,t)$. For $t=1,4,5,6$, $x\ge 10$ and for $t=2,3,7,8,9$, $s=1,2,5,6,7$ respectively, so $y\gg 10$. No solutions exist.
\\
\\
Suppose that $b$ has ten digits; that is, $b=147wxyz247$. The Diophantine equation is $2999w+290x-70y-997z=10000$ and it can be verified from [1] that all solutions produced have at least one variable that is no less than $10$. No solutions exist.
\\
\\
Now suppose that $b$ has eleven digits; that is, $b=147vwxyz247$. The Diophantine equation is $2990w+200x-970y-9997z=100000-29999v$. Going through $v$ from $1$ to $9$ and solving this gives us that only $v=5$ provides one solution with all variables less than $10$. The unique solution is $(v,w,x,y,z)=(5,2,4,7,5)$ so $b=14752475247$.
\\
\\
In particular, we have that for $b$ of any length greater than eleven, $b=1475247\cdots 5247$. But we have arrived at exactly the same situation as in (*). This means that the cycle repeats, and thus $b$ is equal to $147$ followed by blocks of $5247$. Formally, $b=147\cdot 10^{4(k+1)}+5247\sum_{i=0}^k10^{4i}$ and $147$ for $k=0,1,2,\cdots$, and this completes $c_{n-1}=4$.

\paragraph{When $c_{n-1}=6$}

It can be verified that for $n\le 3$, no solutions exist. This can be done by employing linear Diophantine equations. From (6), since $c_n=1$ and $c_1=4$,
\begin{equation*}
-2+3\cdot 10^{n-1}+(3\cdot 10^{n-2}-10)c_{n-1}+\cdots+(30-10^{n-2})c_2+12-4\cdot 10^{n-1}=-10
\end{equation*}
Combining terms and dividing by $10$ gives (for $n>3$)
\begin{equation*}
-10^{n-2}+(3\cdot 10^{n-3}-1)c_{n-1}+\cdots+(3-10^{n-3})c_2=-2 
\end{equation*}
\begin{equation*}
-8-10^{n-2}+(3\cdot 10^{n-3}-1)c_{n-1}+\cdots+(3-10^{n-3})c_2=-10
\end{equation*}
This means that for $c_{n-1}=6$, $-8-6+3c_2$ is divisible by $10$; in other words, $c_2=8$. Hence $b-3=16\cdots\cdots 84\implies b=16\cdots\cdots 87$ and $3b=48\cdots\cdots 61$. In turn, it can be implied that $c_{n-2}=0,1,2,3$ as $16\cdot 3=48$.
\\
\\
If $c_{n-2}=0$, $3b$ ends in $061$ forcing $c_3=6$ since $684\cdot 3=2061$ which is a contradiction as $b$ starts with $160$ and $3b$ starts with $486$. If $c_{n-2}=2$, $3b$ ends in $261$ forcing $c_3=0$ since $84\cdot 3=261$ which is contradictory as $b$ starts with $162$ and $3b$ starts with $480$. Furthermore, if $c_{n-2}=3$, $3b$ ends in $361$ forcing $c_3=7$ since $787\cdot 3=2361$ which is again a contradiction as $b$ starts with $163$ and $3b$ starts with $487<163\cdot 3$. There is no contradiction when $c_{n-2}=1$ because $c_3=3$ works so we can continue down this route.
\\
\\
We are now left with $b=161\cdots\cdots 387\implies b-3=161\cdots\cdots 384\implies 3b=483\cdots\cdots 161$ so again, $c_{n-3}=0,1,2,3$. If $c_{n-3}=0$, $b$ starts with $1610$ so $c_4=0,1,2$ as the fourth digit of $3b$. Now $0387\cdot3=1161$, $1387\cdot3=4161$ and $2387\cdot3=7161$ which is a contradiction as $1,4,7\ne0$ which is what we assumed for $c_{n-3}$ in this case. Similarly, if $c_{n-3}=1$, $b$ starts with $1611$ so $c_4=3,4,5$ as the fourth digit of $3b$. Now $3387\cdot3=10161$, $4387\cdot3=13161$ and $5387\cdot3=16161$ which is a contradiction as $0,3,6\ne1$ which is what we assumed for $c_{n-3}$ in this case. If $c_{n-3}=3$, this forces $c_4=9$ since $1613\cdot 3=4839$ which in turn forces $c_{n-3}=8$ as $9384\cdot 3=28161$. The contradiction is apparent. This means that $b=1612\cdots\cdots 7387$ as it turns out that if $c_{n-3}=2$, $c_4=7$ using the same method as above. ($\dagger$)
\\
\\
Suppose that $b$ has nine digits; that is, $b=1612x7387$ just like in the previous subsection. Then $3b=4836(3x+2)2161$ which means that $20>3x+2\ge 10\implies 3x+2=11,14,17\implies x=3,4,5$. The last digits of $3x+2$ are $1,4,7$ and since $4$ is the only digit such that the last digit of $3x+2$ is the same as $x$, the only solution is $b=161247387$. 
\\
\\
Suppose that $b$ has ten digits; that is, $b=1612xy7387$. From similar reasoning to the above, we deduce that $3(xy)+2=yx+100$. This means that in algebraic terms, $29x-7y=98$. This forces $x=7$ as $29x=7(y+14)$, but this gives $y=15>10$ so no solutions exist.
\\
\\
Suppose that $b$ has eleven digits; that is, $b=1612xyz7387$. The Diophantine equation this time is $299x+20y-97z=998$ and it can be verified from [1] that there are no solutions such that $x,y,z<10$.
\\
\\
Suppose that $b$ has twelve digits; that is, $b=1612wxyz7387$. The Diophantine equation is $2999w+290x-70y-997z=9998$ and similarly it can be verified from [1] that there are no solutions such that $w,x,y,z<10$.
\\
\\
Now suppose that $b$ has thirteen digits; that is, $b=1612vwxyz7387$. The Diophantine equation is $2990w+200x-970y-9997z=99998-29999v$. Going through $v$ from $1$ to $9$ and solving this gives us that only $v=4$ provides one solution with all variables less than $10$. The unique solution is $(v,w,x,y,z)=(4,7,5,2,4)$ so $b=161\,247\,5247\,387$.
\\
\\
In particular, we have that for $b$ of any length greater than thirteen, $b=161\,247\,5247\cdots 387\implies b=161\,247\,52\cdots 387$. But this is exactly where we had arrived at previously in ($\dagger$). This means that the cycle repeats, and thus $b$ is equal to $161247$ followed by blocks of $5247$, followed by $387$ at the end. Formally, $b=161247\cdot 10^{4k+7}+5247\sum_{i=0}^k10^{4i+3}+387$ and $161\,247\,387$ for $k=0,1,2,\cdots$, and this completes $c_{n-1}=6$.

\section{References}

$[1]$ HackMath, (2018). \textit{Integer Diophantine Equations Solver}. Available from:
\\
 \url{https://www.hackmath.net/en/calculator/integer-diophantine-equations-solver}. [Accessed on 3 December 2018].
\\
\\
$[2]$ Nandor, M., (2009). \textit{The On-Line Encyclopedia of Integer Sequences}. OEIS Foundation Inc. Available from: \url{http://oeis.org/A166749}. [Accessed 20 December 2018].
\\
\\
$[3]$ Quora, (2018). \textit{"What are numbers whose sum is the reverse of their product?"} Online posting. Available from: \\\url{https://www.quora.com/What-are-numbers-whose-sum-is-reverse-of-their-product}. [Accessed 20 December 2018].
\\
\\
$[4]$ Lugo, M, (2018). A comment: \textit{When are $a+b$ and $ab$ palindromic for integers $a,b$?} Mathematics Stack Exchange, Stack Exchange Inc. Available from: \\\url{https://math.stackexchange.com/q/2961866/471884} [Accessed 20 December 2018].

\subsection*{Affiliations}

\textsc{Wang Pok Lo, The University of Sheffield, School of Mathematics and Statistics, Hicks Building, Sheffield, S3 7RH, United Kingdom.}

\textit{Email address:} 
\texttt{\href{mailto:wplo1@sheffield.ac.uk}{wplo1@sheffield.ac.uk}}
\\
\\
\textsc{Yuval Paz, The Hebrew University of Jerusalem, Department of Mathematics, E. Safra, Jerusalem, Israel.}

\textit{Email address:} \texttt{\href{mailto:yuval.paz1@mail.huji.ac.il}{yuval.paz1@mail.huji.ac.il}}

\end{document}